\documentclass[11pt]{article}
\usepackage{amssymb,amsthm,amsmath}

\newtheorem{thm}{Theorem}[section]
\newtheorem{prop}[thm]{Proposition}
\newtheorem{lem}[thm]{Lemma}
\newtheorem{cor}[thm]{Corollary}

\theoremstyle{definition}
\newtheorem{df}[thm]{Definition}
\newtheorem{rem}[thm]{Remark}

\newcommand{\N}{\mathbb{N}}
\newcommand{\Z}{\mathbb{Z}}

\newcommand{\R}{\mathbb{R}}
\newcommand{\C}{\mathbb{C}}
\newcommand{\T}{\mathbb{T}}

\newcommand{\Aut}{\operatorname{Aut}}
\newcommand{\Ad}{\operatorname{Ad}}
\newcommand{\id}{\operatorname{id}}
\newcommand{\diag}{\operatorname{diag}}
\newcommand{\Sp}{\operatorname{Sp}}
\newcommand{\Tr}{\operatorname{Tr}}
\newcommand{\Lip}{\operatorname{Lip}}
\newcommand{\Bott}{\operatorname{Bott}}

\newcommand{\ep}{\varepsilon}

\title{Classification of uniformly outer actions of $\Z^2$ \\
on UHF algebras}
\author{Takeshi Katsura 
\thanks{Supported in part by 
the 21st Century COE program at University of Tokyo 
of the Ministry of Education, Culture, Sports, Science and Technology 
of Japan, and by the Fields Institute} \\
Department of Mathematical Sciences \\
University of Tokyo \\
Komaba, Tokyo, 153-8914, Japan 
\and
Hiroki Matui 
\thanks{Supported in part by a grant 
from the Japan Society for the Promotion of Science} \\
Graduate School of Science \\
Chiba University \\
1-33 Yayoi-cho, Inage-ku, Chiba 263-8522, Japan}
\date{}

\begin{document}
\maketitle

\begin{abstract}
We give a complete classification up to cocycle conjugacy 
of uniformly outer actions of $\Z^2$ on UHF algebras. 
In particular, it is shown that 
any two uniformly outer actions of $\Z^2$ 
on a UHF algebra of infinite type are cocycle conjugate. 
We also classify them up to outer conjugacy. 
\end{abstract}

%%%%%%%%%%%%%%%%%%%%%%%%%%%%%%%%%%%%%%%%%%%%%%%%%%%%%%%%%%%%
\section{Introduction}

In \cite{N1}, H. Nakamura introduced 
the Rohlin property for $\Z^N$-actions on unital $C^*$-algebras 
and obtained the Rohlin type theorem 
for $\Z^2$-actions on UHF algebras (see Theorem \ref{RohlinType}). 
He also classified product type actions of $\Z^2$ 
up to cocycle conjugacy. 
In this paper, we extend this result and 
complete classification of uniformly outer actions of $\Z^2$ 
on UHF algebras up to cocycle conjugacy. 

In the theory of operator algebras, 
the problem of classifying group actions has a long history. 
A fundamental tool for the classification theory of actions 
is the (noncommutative) Rohlin property, 
which was first introduced by A. Connes 
for single automorphisms of finite von Neumann algebras (\cite{C}). 
In the framework of $C^*$-algebras, 
A. Kishimoto established a non-commutative Rohlin type theorem 
for single automorphisms of UHF algebras and certain AT algebras, 
and classified them up to outer conjugacy 
(\cite{K1},\cite{K2},\cite{K3}). 
Nakamura used the same technique 
for automorphisms of purely infinite simple nuclear $C^*$-algebras 
and classified them up to outer conjugacy (\cite{N2}). 
Recently, the second named author showed that 
any two outer actions of $\Z^N$ on the Cuntz algebra $\mathcal{O}_2$ 
are cocycle conjugate (\cite{M}). 
In the case of finite group actions, 
M. Izumi introduced a notion of the Rohlin property 
and classified a large class of actions (\cite{I2},\cite{I3}). 
The reader may consult the survey paper \cite{I1} 
for the Rohlin property of automorphisms of $C^*$-algebras. 

The aim of this paper is 
to extend these results to uniformly outer actions of $\Z^2$ 
on UHF algebras. 
More precisely, we will introduce a $K$-theoretical invariant 
for $\Z^2$-actions on a UHF algebra and 
show that it is a complete invariant for cocycle conjugacy. 
When the UHF algebra $A$ is of infinite type, that is, 
$A$ is isomorphic to $A\otimes A$, 
the range of this invariant is zero, and so 
any two actions are cocycle conjugate. 
Note that, in \cite[Theorem 15]{N1}, 
essentially the same invariant was used 
to distinguish infinitely many product type actions 
which are not cocycle conjugate to each other. 

The content of this paper is as follows. 
In Section 2, 
we collect notations and basic facts relevant to this paper. 
Notions of the ultraproduct algebra $A^\infty$ and 
the central sequence algebra $A_\infty$ will help our analysis. 
For $\Z^2$-actions on unital $C^*$-algebras, 
we recall the definition of the Rohlin property from \cite{N1}. 
In Section 3, we introduce the notion of admissible cocycles. 
For an almost cocycle of an approximately inner $\Z^2$-action 
on an AF algebra, 
we can associate an element of the $K_0$-group. 
An almost cocycle is said to be admissible, 
if the associated element in the $K_0$-group is zero. 
Section 4 is devoted to the cohomology vanishing theorem 
for $\Z^2$-actions on UHF algebras with the Rohlin property. 
One of the difficulties 
in the study of $\Z^2$-actions on $C^*$-algebras is that 
one has to formulate a two-dimensional version of Berg's technique 
in order to obtain the so-called cohomology vanishing theorem. 
In our situation, 
we need to construct a certain homotopy of unitaries, 
and admissibility of a cocycle is necessary for that. 
In Section 5, 
a $K$-theoretical invariant for $\Z^2$-actions is introduced. 
We also observe that 
the invariant is exhausted by uniformly outer actions. 
In Section 6, 
for given two uniformly outer $\Z^2$-actions $\alpha,\beta$ 
with the same invariant, 
we show the existence of cocycles 
which transform $\alpha$ to $\beta$ approximately. 
Kishimoto's classification result of single automorphisms 
is used to obtain the admissibility of cocycles. 
Finally, the main theorems are given. 
D. E. Evans and A. Kishimoto introduced in \cite{EK} 
an intertwining argument for automorphisms, 
which is an equivariant version of Elliott's intertwining argument 
for classification of $C^*$-algebras. 
By using the Evans-Kishimoto intertwining argument, 
we get the classification theorem up to cocycle conjugacy. 
The classification up to outer conjugacy is also obtained.

%%%%%%%%%%%%%%%%%%%%%%%%%%%%%%%%%%%%%%%%%%%%%%%%%%%%%%%%%%%%
\section{Preliminaries}

Let $A$ be a unital $C^*$-algebra. 
We denote by $U(A)$ the group of unitaries in $A$. 
For $u\in U(A)$, we let $\Ad u(a)=uau^*$ for $a\in A$ and 
call it an inner automorphism of $A$. 
For any $a,b\in A$, we write $[a,b]=ab-ba$ 
and call it the commutator of $a$ and $b$. 
We let $\log$ be the standard branch 
defined on the complement of the negative real axis. 

First, we would like to recall the classification result of 
single automorphisms of UHF algebras. 

\begin{df}[{\cite[Definition 1.2]{K1}}]
An automorphism $\alpha$ of a unital $C^*$-algebra $A$ is 
said to be uniformly outer 
if for any $a\in A$, any non-zero projection $p\in A$ and any $\ep>0$, 
there exist projections $p_1,p_2,\dots,p_n$ in $A$ such that 
$p=\sum p_i$ and $\lVert p_ia\alpha(p_i)\rVert<\ep$ 
for all $i=1,2,\dots,n$. 
\end{df}

We say that 
an action $\alpha$ of a discrete group on $A$ is uniformly outer 
if $\alpha_n$ is uniformly outer for every element $n$ of the group 
other than the identity element. 
The following result can be found in \cite{K1} 
(see also \cite{BEK}). 

\begin{thm}[{\cite[Theorem 1.3]{K1}}]\label{cc}
Let $A$ be a UHF algebra and 
let $\alpha,\beta\in\Aut(A)$ be automorphisms with 
which the associated actions of $\Z$ are uniformly outer. 
Then for any $\ep>0$, 
there exist an automorphism $\sigma$ and a unitary $u\in A$ 
such that 
\[ \lVert u-1\rVert<\ep \ \text{ and } \ 
\Ad u\circ\alpha=\sigma\circ\beta\circ\sigma^{-1}. \]
\end{thm}

We remark that 
Kishimoto generalized the result above to AT algebras 
in \cite{K2} and \cite{K3}. 
In this paper, 
we would like to study a $\Z^2$ version of this classification theorem 
on UHF algebras. 
A key point for the proof of the theorem above is the Rohlin property 
of automorphisms or $\Z$-actions. 
We would like to recall 
the definition of the Rohlin property for $\Z^2$-actions 
on unital $C^*$-algebras (see \cite[Section 2]{N1} or \cite[Section 3]{I1}). 

For $m=(m_1,m_2)$ and $n=(n_1,n_2)$ in $\Z^2$, 
$m\leq n$ means $m_1\leq n_1$ and $m_2\leq n_2$. 
For $m=(m_1,m_2)\in\N^2$, we let 
\[ m\Z^2=\{(m_1n_1,m_2n_2)\in\Z^2
\mid (n_1,n_2)\in\Z^2\}. \]
For simplicity, we denote $\Z^2/m\Z^2$ by $\Z_m$. 
Moreover, we may identify $\Z_m=\Z^2/m\Z^2$ with 
\[ \{(n_1,n_2)\in\Z^2
\mid 0\leq n_i\leq m_i-1\text{ for each }i=1,2\}. \]
The canonical basis of $\Z^2$ as well as its image in $\Z_m$ 
is denoted by $\xi_1=(1,0)$ and $\xi_2=(0,1)$. 

\begin{df}\label{Rohlin}
Let $\alpha$ be an action of $\Z^2$ on a unital $C^*$-algebra $A$. 
Then $\alpha$ is said to have the Rohlin property, 
if for any $m\in\N$ there exist $R\in\N$ and 
$m^{(1)},m^{(2)},\dots,m^{(R)}\in\N^2$ 
with $m^{(1)},\dots,m^{(R)}\geq(m,m)$ 
satisfying the following: 
For any finite subset $\mathcal{F}$ of $A$ and $\ep>0$, 
there exists a family of projections 
\[ e^{(r)}_g \qquad (r=1,2,\dots,R, \ g\in\Z_{m^{(r)}}) \]
in $A$ such that 
\[ \sum_{r=1}^R\sum_{g\in\Z_{m^{(r)}}}e^{(r)}_g=1, 
\quad \lVert[a,e^{(r)}_g]\rVert<\ep, 
\quad \lVert\alpha_{\xi_i}(e^{(r)}_g)-e^{(r)}_{g+\xi_i}\rVert<\ep \]
for any $a\in\mathcal{F}$, $r=1,2,\dots,R$, $i=1,2$ 
and $g\in\Z_{m^{(r)}}$, 
where $g+\xi_i$ is understood modulo $m^{(r)}\Z^2$. 
\end{df}

Nakamura proved the following Rohlin type theorem 
for $\Z^2$-actions on UHF algebras in \cite{N1}. 

\begin{thm}[{\cite[Theorem 3]{N1}}]\label{RohlinType}
Let $\alpha$ be an action of $\Z^2$ on a UHF algebra $A$. 
Then the following conditions are equivalent. 
\begin{enumerate}
\item $\alpha$ has the Rohlin property. 
\item $\alpha$ is uniformly outer. 
\end{enumerate}
\end{thm}

In this paper, we deal with central sequence algebras, 
which simplify the arguments a little. 
Let $A$ be a separable $C^*$-algebra. 
We set 
\[ c_0(A)=\{(a_n)\in\ell^\infty(\N,A)\mid
\lim_{n\to\infty}\lVert a_n\rVert=0\}, \]
\[ A^\infty=\ell^\infty(\N,A)/c_0(A). \]
We identify $A$ with the $C^*$-subalgebra of $A^\infty$ 
consisting of equivalence classes of constant sequences. 
We let 
\[ A_\infty=A^\infty\cap A'. \]
When $\alpha$ is an automorphism of $A$ or 
an action of a discrete group on $A$, 
we can consider its natural extension on $A^\infty$ and $A_\infty$. 
We denote it by the same symbol $\alpha$. 

Using these notations, 
we can restate the definition of the Rohlin property 
of a $\Z^2$-action $\alpha$ on a unital separable $C^*$-algebra $A$ 
as follows. 
For any $m\in\N$ there exist $R\in\N$, 
$m^{(1)},m^{(2)},\dots,m^{(R)}\in\N^2$ 
with $m^{(1)},\dots,m^{(R)}\geq(m,m)$ and 
a family of projections 
\[ e^{(r)}_g \qquad (r=1,2,\dots,R, \ g\in\Z_{m^{(r)}}) \]
in $A_\infty=A^\infty\cap A'$ such that 
\[ \sum_{r=1}^R\sum_{g\in\Z_{m^{(r)}}}e^{(r)}_g=1, 
\quad \alpha_{\xi_i}(e^{(r)}_g)=e^{(r)}_{g+\xi_i} \]
for any $r=1,2,\dots,R$, $i=1,2$ and $g\in\Z_{m^{(r)}}$, 
where $g+\xi_i$ is understood modulo $m^{(r)}\Z^2$. 

\bigskip

Let $\xi_1=(1,0)$ and $\xi_2=(0,1)$ be the canonical basis of $\Z^2$. 
Let $\alpha$ be an action of $\Z^2$ on a unital $C^*$-algebra $A$. 
A family of unitaries $\{u_n\}_{n\in\Z^2}$ in $A$ is called 
an $\alpha$-cocycle, 
if 
\[ u_n\alpha_n(u_m)=u_{n+m} \]
for all $n,m\in\Z^2$. 
If a pair of unitaries $u_1,u_2\in U(A)$ satisfies 
\[ u_1\alpha_{\xi_1}(u_2)=u_2\alpha_{\xi_2}(u_1), \]
then it determines uniquely an $\alpha$-cocycle $\{u_n\}_{n\in\Z^2}$ 
such that $u_{\xi_i}=u_i$. 
We may also call the pair $\{u_1,u_2\}$ an $\alpha$-cocycle. 
An $\alpha$-cocycle $\{u_n\}_n$ in $A$ is called a coboundary, 
if there exists $v\in U(A)$ such that 
\[ u_n=v\alpha_n(v^*) \]
for all $n\in\Z^2$, or equivalently, 
if 
\[ u_{\xi_i}=v\alpha_{\xi_i}(v^*) \]
for each $i=1,2$. 

When $\{u_n\}_{n\in\Z^2}$ is an $\alpha$-cocycle, 
it turns out that 
a new action $\tilde{\alpha}$ of $\Z^2$ on $A$ can be defined by 
\[ \tilde{\alpha}_n(x)=\Ad u_n\circ\alpha_n(x)=u_n\alpha_n(x)u_n^* \]
for each $x\in A$. 
We call $\tilde{\alpha}$ 
the perturbed action of $\alpha$ by $\{u_n\}_n$. 

Two actions $\alpha$ and $\beta$ of $\Z^2$ on $A$ 
are said to be cocycle conjugate, 
if there exists an $\alpha$-cocycle $\{u_n\}_n$ in $A$ such that 
the perturbed action of $\alpha$ by $\{u_n\}_n$ is conjugate to $\beta$. 
The purpose of this paper is 
to classify uniformly outer actions of $\Z^2$ on UHF algebras 
up to cocycle conjugacy.

%%%%%%%%%%%%%%%%%%%%%%%%%%%%%%%%%%%%%%%%%%%%%%%%%%%%%%%%%%%%
\section{Admissible cocycles}

Throughout this section, 
we let $A$ be a unital AF algebra and 
let $\alpha$ be an action of $\Z^2$ on $A$ 
such that $\alpha_n$ is approximately inner for every $n\in\Z^2$. 
To simplify notation, 
we may denote $\alpha_{\xi_i}$ by $\alpha_i$ for each $i=1,2$. 

A pair of unitaries $\{u_1,u_2\}$ in $A$ is called 
an almost $\alpha$-cocycle, when 
\[ \lVert u_1\alpha_1(u_2)-u_2\alpha_2(u_1)\rVert<1. \]
For such a pair $u_1,u_2$, 
we can construct a closed path of unitaries as follows. 
We put 
\[ a=\frac{1}{2\pi\sqrt{-1}}\log(u_1\alpha_1(u_2)(u_2\alpha_2(u_1))^*), \]
and set a path $k:[0,1]\to U(A)$ 
from $u_2\alpha_2(u_1)$ to $u_1\alpha_1(u_2)$ by 
\[ k(t)=e^{2\pi\sqrt{-1}ta}u_2\alpha_2(u_1). \]
Let $h_i:[0,1]\to U(A)$ be a path from $1$ to $u_i$ for each $i=1,2$. 
Then 
\[ \widetilde{h}_1:t\mapsto u_2\alpha_2(h_1(t)) \]
is a path from $u_2$ to $u_2\alpha_2(u_1)$ and 
\[ \widetilde{h}_2:t\mapsto u_1\alpha_1(h_2(t)) \]
is a path from $u_1$ to $u_1\alpha_1(u_2)$. 
By connecting these paths, 
we obtain a closed path $H:[0,1]\to U(A)$. 
More precisely, 
\[ H(t)=\begin{cases}
h_1(5t) & \text{ if }0\leq t\leq1/5 \\
\widetilde{h}_2(5t-1) & \text{ if }1/5\leq t\leq2/5 \\
k(3-5t)  & \text{ if }2/5\leq t\leq3/5 \\
\widetilde{h}_1(4-5t) & \text{ if }3/5\leq t\leq4/5 \\
h_2(5-5t) & \text{ if }4/5\leq t\leq1. \end{cases} \]
Since $\alpha$ is approximately inner, 
the $K_1$-class of this closed path $H$ in $K_1(SA)\cong K_0(A)$ 
does not depend on the choice of $h_1$ and $h_2$. 
We denote this element in $K_0(A)$ by 
\[ \kappa(u_1,u_2,\alpha_1,\alpha_2)\in K_0(A). \]
Clearly $\kappa(u_1,u_2,\alpha_1,\alpha_2)$ is homotopy invariant 
in the following sense: 
if $u_i:[0,1]\to U(A)$ is a path of unitaries satisfying 
\[ \lVert u_1(t)\alpha_1(u_2(t))
-u_2(t)\alpha_2(u_1(t))\rVert<1 \]
for every $t\in[0,1]$, then 
\[ \kappa(u_1(0),u_2(0),\alpha_1,\alpha_2)
=\kappa(u_1(1),u_2(1),\alpha_1,\alpha_2). \]

\begin{figure}
\begin{center}
\begin{picture}(400,170)

\put(130,30){\line(1,0){140}}
\put(130,60){\line(0,1){100}}
\put(130,160){\line(1,0){140}}
\put(270,160){\line(0,-1){130}}

\multiput(130,30)(140,0){2}{\circle*{5}}
\put(130,60){\circle*{5}}
\multiput(130,160)(140,0){2}{\circle*{5}}

\qbezier[10](130,30)(130,40)(130,60)

\put(275,160){$1$}
\put(196,164){$h_1$}
\put(114,160){$u_1$}
\put(84,58){$u_1\alpha_1(u_2)$}
\put(100,18){$u_2\alpha_2(u_1)$}
\put(270,20){$u_2$}
\put(275,90){$h_2$}

\end{picture}
\caption{The definition of $\kappa(u_1,u_2,\alpha_1,\alpha_2)$}
\end{center}
\end{figure}
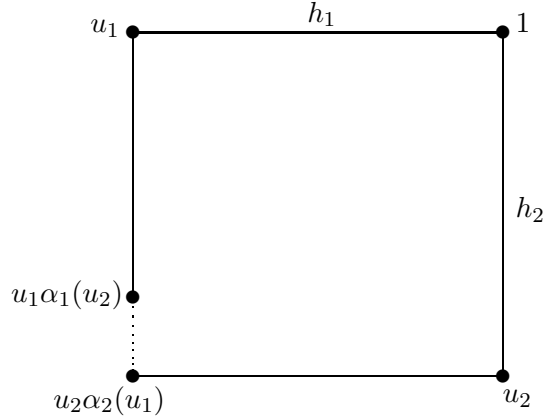

\begin{lem}\label{tracevalue}
Let $\{u_1,u_2\}$ be an almost $\alpha$-cocycle with 
$\lVert u_1\alpha_1(u_2)-u_2\alpha_2(u_1)\rVert<\ep<1$. 
Then, 
\[ \lvert\tau(\kappa(u_1,u_2,\alpha_1,\alpha_2))\rvert
<\frac{\arcsin\ep}{2\pi} \]
for any tracial state $\tau\in T(A)$. 
\end{lem}
\begin{proof}
Let $k$, $h_1$, $h_2$, $\widetilde{h}_1$, $\widetilde{h}_2$ and $H$ 
be the paths of unitaries as above. 
We may assume that $H$ is piecewise smooth. 
For any tracial state $\tau\in T(A)$, we have 
\[ \tau(\kappa(u_1,u_2,\alpha_1,\alpha_2))
=\frac{1}{2\pi\sqrt{-1}}\int_0^1\tau(\dot{H}(t)H(t)^*)\,dt. \]
Since $\alpha_2$ is approximately inner, 
$\tau\circ\alpha_2$ is equal to $\tau$. 
Therefore the contribution from $h_1$ and $\widetilde{h}_1$ cancels out. 
Similarly, 
the contribution from $h_2$ and $\widetilde{h}_2$ also cancels out. 
Therefore we get $\tau(\kappa(u_1,u_2,\alpha_1,\alpha_2))=\tau(a)$, 
where 
\[ a=\frac{1}{2\pi\sqrt{-1}}\log(u_1\alpha_1(u_2)(u_2\alpha_2(u_1))^*). \]
Since $\lVert a\rVert<\pi^{-1}\arcsin(\ep/2)<(2\pi)^{-1}\arcsin\ep$, 
the proof is completed. 
\end{proof}

Suppose that 
a family of unitaries $\{u_n\}_{n\in\Z^2}$ in $A^\infty$ 
is an $\alpha$-cocycle.  
Let $(u_{\xi_i}(k))_{k\in\N}$ be 
a representing sequence in $\ell^\infty(\N,A)$ of $u_{\xi_i}$. 
We may assume that $u_{\xi_i}(k)$ is in $U(A)$. 
For sufficiently large $k\in\N$, 
the pair $\{u_{\xi_1}(k),u_{\xi_2}(k)\}$ is an almost $\alpha$-cocycle in $A$. 
In addition, from the lemma above, 
we can see that 
\[ -1\leq\kappa(u_{\xi_1}(k),u_{\xi_2}(k),\alpha_{\xi_1},\alpha_{\xi_2})
\leq1 \]
for sufficiently large $k\in\N$. 
Hence 
\[ \kappa^\infty(u_{\xi_1},u_{\xi_2},\alpha_{\xi_1},\alpha_{\xi_2})
=(\kappa(u_{\xi_1}(k),u_{\xi_2}(k),\alpha_{\xi_1},\alpha_{\xi_2}))_{k\in\N}
\in K_0(A^\infty) \]
is well-defined. 

\begin{df}
Let $A$ be a unital AF algebra and 
let $\alpha$ be an action of $\Z^2$ on $A$ 
such that $\alpha_n$ is approximately inner for every $n\in\Z^2$. 
\begin{enumerate}
\item An $\alpha$-cocycle $\{u_n\}_{n\in\Z^2}$ in $A$ is 
said to be admissible, 
if $\kappa(u_{\xi_1},u_{\xi_2},\alpha_{\xi_1},\alpha_{\xi_2})$ 
is zero in $K_0(A)$. 
\item An $\alpha$-cocycle $\{u_n\}_{n\in\Z^2}$ in $A^\infty$ is 
said to be admissible, 
if $\kappa^\infty(u_{\xi_1},u_{\xi_2},\alpha_{\xi_1},\alpha_{\xi_2})$ 
is zero in $K_0(A^\infty)$. 
\end{enumerate}
\end{df}
Observe that if $\alpha$-cocycle is a coboundary, 
then it is admissible. 
When $K_0(A)$ has no non-trivial infinitesimal elements, 
from Lemma \ref{tracevalue}, 
we can see that any $\alpha$-cocycle in $A$ is admissible. 

Let $\{u_1,u_2\}$ be an almost $\alpha$-cocycle in $A$. 
For $j\in\N$ and $i=1,2$, we define $u^{(j)}_i\in U(A)$ by 
\[ u^{(1)}_i=u_i\text{ and }
u^{(j+1)}_i=u^{(j)}_i\alpha_i^j(u_i). \]
Suppose that 
\[ \lVert u^{(j)}_1\alpha_1^j(u^{(k)}_2)
-u^{(k)}_2\alpha_2^k(u^{(j)}_1)\rVert<1 \]
for every $j=1,2,\dots,m$ and $k=1,2,\dots,n$. 
Thus, the pair $\{u^{(j)}_1,u^{(k)}_2\}$ is an almost cocycle 
for the $\Z^2$-action generated by $\alpha_1^j$ and $\alpha_2^k$. 

\begin{lem}\label{power}
In the setting above, one has 
\[ \kappa(u^{(m)}_1,u^{(n)}_2,\alpha_1^m,\alpha_2^n)
=nm\kappa(u_1,u_2,\alpha_1,\alpha_2) \]
in $K_0(A)$. 
\end{lem}
\begin{proof}
It suffices to show 
\[ \kappa(u^{(j+1)}_1,u^{(k)}_2,\alpha_1^{j+1},\alpha_2^k)
=\kappa(u^{(j)}_1,u^{(k)}_2,\alpha_1^j,\alpha_2^k)
+\kappa(u_1,u^{(k)}_2,\alpha_1,\alpha_2^k). \]
Let $h_1:[0,1]\to U(A)$ be a path from $1$ to $u_1$ and 
let $h_2:[0,1]\to U(A)$ be a path from $1$ to $u^{(k)}_2$. 
Let $H$ be a closed path of unitaries obtained by connecting 
a short path from $u^{(j)}_1\alpha_1^j(u^{(k)}_2)$ 
to $u^{(k)}_2\alpha_2^k(u^{(j)}_1)$, 
a short path from $u^{(j+1)}_1\alpha_1^{j+1}(u^{(k)}_2)$ 
to $u^{(k)}_2\alpha_2^k(u^{(j+1)}_1)$ and 
the following four paths: 
\[ t\mapsto u^{(j)}_1\alpha_1^j(h_1(t)), 
\qquad t\mapsto u^{(j)}_1\alpha_1^j(h_2(t)) \]
\[ t\mapsto u^{(k)}_2\alpha_2^k(u^{(j)}_1\alpha_1^j(h_1(t))), 
\qquad t\mapsto u^{(j+1)}_1\alpha_1^{j+1}(h_2(t)). \]
Then, the $K_0$-class of this closed path $H$ 
in $K_1(SA)\cong K_0(A)$ is equal to 
\[ \kappa(u^{(j+1)}_1,u^{(k)}_2,\alpha_1^{j+1},\alpha_2^k)
-\kappa(u^{(j)}_1,u^{(k)}_2,\alpha_1^j,\alpha_2^k). \]
Let $H'$ be a closed path of unitaries obtained by connecting 
a short path 
from $u^{(k)}_2\alpha_2^k(u_1)$ to $u_1\alpha_1(u^{(k)}_2)$ and 
the following four paths: 
\[ t\mapsto h_1(t), 
\qquad t\mapsto h_2(t) \]
\[ t\mapsto u^{(k)}_2\alpha_2^k(h_1(t)), 
\qquad t\mapsto u_1\alpha_1(h_2(t)). \]
By definition, the $K_0$-class of $H'$ equals 
$\kappa(u_1,u^{(k)}_2,\alpha_1,\alpha_2^k)$. 
It follows from 
\[ u^{(j)}_1\alpha_1^j(u^{(k)}_2)
(u^{(k)}_2\alpha_2^k(u^{(j)}_1))^*
u^{(k)}_2\alpha_2^k(u^{(j+1)}_1)
=u^{(j)}_1\alpha_1^j(u^{(k)}_2\alpha_2^k(u_1)) \]
that $H$ is homotopic 
to the closed path $u^{(j)}_1\alpha_1^j(H'(\cdot))$, 
which implies the desired equality. 
\end{proof}

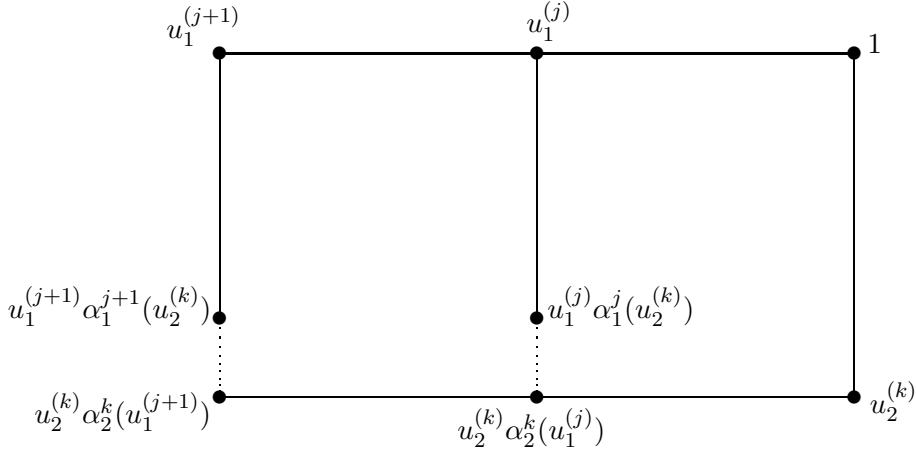
\begin{figure}
\begin{center}
\begin{picture}(400,170)

\put(80,30){\line(1,0){240}}
\put(80,60){\line(0,1){100}}
\put(80,160){\line(1,0){240}}
\put(200,160){\line(0,-1){100}}
\put(320,160){\line(0,-1){130}}

\multiput(80,160)(120,0){3}{\circle*{5}}
\multiput(80,60)(120,0){2}{\circle*{5}}
\multiput(80,30)(120,0){3}{\circle*{5}}

\qbezier[10](80,30)(80,40)(80,60)
\qbezier[10](200,30)(200,40)(200,60)

\put(325,160){$1$}
\put(196,168){$u^{(j)}_1$}
\put(60,166){$u^{(j+1)}_1$}
\put(0,60){$u^{(j+1)}_1\alpha_1^{j+1}(u^{(k)}_2)$}
\put(10,20){$u^{(k)}_2\alpha_2^k(u^{(j+1)}_1)$}
\put(204,60){$u^{(j)}_1\alpha_1^j(u^{(k)}_2)$}
\put(170,14){$u^{(k)}_2\alpha_2^k(u^{(j)}_1)$}
\put(326,24){$u^{(k)}_2$}

\end{picture}
\caption{$\kappa(u^{(j+1)}_1,u^{(k)}_2,\alpha_1^{j+1},\alpha_2^k)$ 
and $\kappa(u^{(j)}_1,u^{(k)}_2,\alpha_1^j,\alpha_2^k)$}
\end{center}
\end{figure}

Let $\{u_1,u_2\}$ be an almost $\alpha$-cocycle in $A$ such that 
\[ \lVert u_1\alpha_1(u_2)-u_2\alpha_2(u_1)\rVert<\ep. \]
Suppose that a unitary $v$ in $A$ satisfies 
\[ \lVert v-u_1\alpha_1(v)u_1^*\rVert<\delta. \]
Then we have 
\[ \lVert u_1\alpha_1(vu_2)-vu_2\alpha_2(u_1)\rVert<\ep+\delta. \]
Thus, if $\delta$ is sufficiently small, then 
the pair $\{u_1,vu_2\}$ is again an almost $\alpha$-cocycle. 

\begin{lem}\label{adjust}
In the setting above, we have 
\[ \kappa(u_1,vu_2,\alpha_1,\alpha_2)
=\kappa(u_1,u_2,\alpha_1,\alpha_2)
+\kappa(1,v,\Ad u_1\circ\alpha_1,\id) \]
in $K_0(A)$. 
\end{lem}
\begin{proof}
One can prove this statement in a similar fashion to Lemma \ref{power}. 
\end{proof}

The following two lemmas are key ingredients 
for the cohomology vanishing theorem in the next section. 
We denote the $\ell^\infty$-norm on $\R^2$ by $\lVert\cdot\rVert$. 
We put 
\[ E=\{t\in\R^2\mid\lVert t\rVert\leq1\} \]
and 
\[ \partial E=\{t\in\R^2\mid\lVert t\rVert=1\}. \]
For a Lipschitz continuous function $f$, 
we denote its Lipschitz constant by $\Lip(f)$. 

\begin{lem}\label{zONboundary1}
Let $\ep>0$ be a sufficiently small number. 
Suppose that a pair of unitaries $u_1,u_2\in A$ satisfies 
\[ \lVert u_1\alpha_1(u_2)-u_2\alpha_2(u_1)\rVert<\ep \]
and $\kappa(u_1,u_2,\alpha_1,\alpha_2)=0$. 
Then there exists a continuous map $z:\partial E\to U(A)$ 
such that the following hold. 
\begin{enumerate}
\item $z(1,1)=1$. 
\item $\lVert z(-1,t)-u_1\alpha_1(z(1,t))\rVert<\ep$ and 
$\lVert z(t,-1)-u_2\alpha_2(z(t,1))\rVert<\ep$ for all $t\in[-1,1]$. 
\item $\Lip(z)\leq4$. 
\item The $K_1$-class of $z$ in $K_1(C(\partial E)\otimes A)$ is zero. 
\end{enumerate}
\end{lem}
\begin{proof}
This is immediate from 
the definition of $\kappa(u_1,u_2,\alpha_1,\alpha_2)$ 
and the fact that 
any unitary in an AF algebra can be connected to the identity 
by a path of unitaries with length less than $\pi+\delta$, 
where $\delta$ is an arbitrary positive real number. 
We leave the details to the reader. 
\end{proof}

\begin{lem}\label{zONboundary2}
Let $A_0\subset A$ be a unital finite dimensional $C^*$-subalgebra. 
For any $\ep>0$, there exist $\delta>0$ and 
a unital finite dimensional $C^*$-subalgebra $A_1\subset A$ 
such that the following holds. 
Suppose that a pair of unitaries $u_1,u_2\in A\cap A_1'$ satisfies 
\[ \lVert u_1\alpha_1(u_2)-u_2\alpha_2(u_1)\rVert<\delta \]
and $\kappa(u_1,u_2,\alpha_1,\alpha_2)=0$. 
Then there exists a continuous map $z:\partial E\to U(A\cap A_0')$ 
such that the following hold. 
\begin{enumerate}
\item $z(1,1)=1$. 
\item $\lVert z(-1,t)-u_1\alpha_1(z(1,t))\rVert<\ep$ and 
$\lVert z(t,-1)-u_2\alpha_2(z(t,1))\rVert<\ep$ for all $t\in[-1,1]$. 
\item $\Lip(z)\leq4$. 
\item The $K_1$-class of $z$ in $K_1(C(\partial E)\otimes A)$ is zero. 
\end{enumerate}
\end{lem}
\begin{proof}
Choose sufficiently small $\delta>0$. 
Take a unital finite dimensional $C^*$-subalgebra $A_1$ so that 
for any unitary $w\in A_0$, there exist $a,b,c\in A_1$ such that 
\[ \lVert w-a\rVert<\delta, \quad 
\lVert \alpha_1^{-1}(w)-b\rVert<\delta \ \text{ and } \ 
\lVert \alpha_2^{-1}(w)-c\rVert<\delta. \]
Suppose that a pair of unitaries $u_1,u_2\in A\cap A_1'$ satisfies 
\[ \lVert u_1\alpha_1(u_2)-u_2\alpha_2(u_1)\rVert<\delta \]
and $\kappa(u_1,u_2,\alpha_1,\alpha_2)=0$. 
Since $A\cap A_1'$ is also an AF algebra, for each $i=1,2$, 
there exists a path of unitaries $h_i$ in $A\cap A_1'$ 
from $1$ to $u_i$ with length less than $\pi+\delta$. 
By using these paths, in the same way as the lemma above, 
we can find a continuous map $z_0:\partial E\to U(A)$ 
such that the following hold. 
\begin{itemize}
\item $z_0(1,1)=1$. 
\item For any $w\in U(A_0)$ and $(s,t)\in\partial E$, 
$\lVert[w,z_0(s,t)]\rVert<6\delta$. 
\item $\lVert z_0(-1,t)-u_1\alpha_1(z_0(1,t))\rVert<\delta$ and 
$\lVert z_0(t,-1)-u_2\alpha_2(z_0(t,1))\rVert<\delta$ 
for all $t\in[-1,1]$. 
\item $\Lip(z_0)\leq\pi+2\delta$. 
\item The $K_1$-class of $z_0$ in $K_1(C(\partial E)\otimes A)$ is zero. 
\end{itemize}
Using the Haar measure on the compact group $U(A_0)$, 
we define $z_1:\partial E\to A\cap A_0'$ by 
\[ z_1(s,t)=\int wz_0(s,t)w^*\,dw \]
for every $(s,t)\in\partial E$. 
Then one has $\lVert z_1(s,t)-z_0(s,t)\rVert<6\delta$ 
and $\Lip(z_1)\leq\pi+2\delta$. 
Let 
\[ z_1(s,t)=z(s,t)\lvert z_1(s,t)\rvert \]
be the polar decomposition of $z_1(s,t)$ in $A\cap A_0'$. 
It is not so hard to see that $z(s,t)$ is the desired map. 
\end{proof}

%%%%%%%%%%%%%%%%%%%%%%%%%%%%%%%%%%%%%%%%%%%%%%%%%%%%%%%%%%%%
\section{Cohomology vanishing}

We let $M_n$ denote the $C^*$-algebra of all $n$ by $n$ matrices, 
and $\T$ denote the unit circle of $\C$. 
For a Lipschitz continuous function $f$, 
its Lipschitz constant is denoted by $\Lip(f)$. 

\begin{lem}\label{Lipeigenvalue}
Let $u:[0,1]\to M_n$ be a path of unitaries such that $\Lip(u)\leq C$. 
Suppose that a continuous function $\lambda:[0,1]\to\T$ satisfies 
$\lambda(t)\in\Sp(u(t))$ for any $t\in[0,1]$. 
Then, we have $\Lip(\lambda)\leq C$. 
\end{lem}
\begin{proof}
It suffices to show that 
for each $t\in[0,1]$, there exists $\delta>0$ such that 
$\lvert\lambda(s)-\lambda(t)\rvert\leq C\lvert s-t\rvert$ 
holds for all $s\in(t-\delta,t+\delta)$. 
Take $t\in[0,1]$. 
There exists $\delta>0$ such that 
\[ \lvert \lambda(t)-\mu\rvert>2C\delta \]
for any $\mu\in\Sp(u(t))\setminus\{\lambda(t)\}$. 
From $\Lip(u)\leq C$, 
one has $\lVert u(s)-u(t)\rVert\leq C\lvert s-t\rvert$. 
Therefore, 
the distance from $\lambda(s)$ to $\Sp(u(t))$ is 
not more than $C\lvert s-t\rvert$. 
Since $\lambda$ is continuous and $C\lvert s-t\rvert<C\delta$ 
for $s\in(t-\delta,t+\delta)$, 
we get $\lvert \lambda(s)-\lambda(t)\rvert\leq C\lvert s-t\rvert$ 
for all $s\in(t-\delta,t+\delta)$. 
\end{proof}

The following lemma is 
an immediate consequence of \cite[Corollary 6.15]{ELP} 
(and is also a special case of \cite[Lemma 3.8]{L}). 
When $v$ and $w$ are almost commuting unitaries in $M_n$, 
we denote the Bott element associated with $v$ and $w$ 
by $\Bott(v,w)\in K_0(M_n)\cong\Z$. 
The Bott element $\Bott(v,w)$ can be calculated as follows. 
We put 
\[ a=\frac{1}{2\pi\sqrt{-1}}\log(vwv^*w^*). \]
Since $\det(vwv^*w^*)=1$, we have $\Tr(a)\in\Z$, 
where $\Tr$ is the unnormalized trace on $M_n$. 
Then $\Bott(v,w)$ is equal to $\Tr(a)$. 

\begin{lem}\label{SuperHomotopy}
For any $\ep>0$, there exists $\delta>0$ 
such that the following holds. 
For any $n\in\N$ and unitaries $v,w\in M_n$ satisfying 
\[ \lVert[v,w]\rVert<\delta \text{ and } \Bott(v,w)=0, \]
there exists a continuous path of unitaries 
$w:[0,1]\to M_n$ such that 
\[ w(0)=1, \quad w(1)=w, \quad \Lip(w)\leq\pi+\ep \]
and 
\[ \lVert[v,w(t)]\rVert<\ep \]
for all $t\in[0,1]$. 
\end{lem}

\begin{prop}\label{Lipshrink1}
For any $C>0$ and $\ep>0$, 
there exists $C'>0$ such that the following holds. 
Let $n\in\N$ and 
let $u:[0,1]\to M_n$ be a path of unitaries such that 
$u(0)=u(1)=1$ and $\Lip(u)\leq C$. 
If the $K_1$-class of $u$ in $K_1(C_0(0,1)\otimes M_n)$ is zero, 
then we can find a path of self-adjoint elements 
$h:[0,1]\to M_n$ such that the following are satisfied. 
\begin{enumerate}
\item $\lVert u(t)-e^{2\pi\sqrt{-1}h(t)}\rVert<\ep$ 
for all $t\in[0,1]$. 
\item $h(0)=h(1)=0$. 
\item $\Lip(h)\leq C'$. 
\end{enumerate}
\end{prop}
\begin{proof}
By applying Lemma \ref{SuperHomotopy} to $\ep/2$, 
we get $\delta>0$. We may assume that $\delta$ is less than $\ep/2$. 
Choose $L\in\N$ so that $2C/L<\delta$. 
Put $C'=2CL/3+C/6$. 

Suppose that $u:[0,1]\to M_n$ satisfies the conditions in the statement. 
We would like to construct a path of self-adjoint elements $h:[0,1]\to M_n$ 
satisfying (1), (2) and (3). 
One can find continuous functions $\lambda_1,\lambda_2,\dots,\lambda_n$ 
from $[0,1]$ to $\T$ such that 
$u(t)$ is unitarily equivalent to 
$\diag(\lambda_1(t),\lambda_2(t),\dots,\lambda_n(t))$ for each $t\in[0,1]$. 
Since the $K_1$-class of $u$ in $K_1(C_0(0,1)\otimes M_n)$ is zero, 
we may assume that the rotation number of $\lambda_i$ is zero 
for every $i=1,2,\dots,n$. 
Thus, there exist continuous functions $g_1,g_2,\dots,g_n$ 
from $[0,1]$ to $\R$ such that 
$g_i(0)=g_i(1)=0$ and $\lambda_i(t)=e^{2\pi\sqrt{-1}g_i(t)}$ 
for all $i=1,2,\dots,n$ and $t\in[0,1]$. 
By Lemma \ref{Lipeigenvalue}, we have $\Lip(\lambda_i)\leq C$. 
It follows from $6<2\pi$ that 
the Lipschitz constant of $g_i$ is less than $C/6$. 
In particular, 
\[ \max\{\lvert g_i(t)\rvert\mid t\in[0,1], \ i=1,2,\dots,n\}
\leq\frac{C}{12}. \]

We divide $[0,1]$ to closed intervals 
$I_k=[k/L,(k+1)/L]$ for $k=0,1,\dots,L-1$ 
and construct the path $h$ on each interval. 
For each $k=0,1,\dots,L$, we set $u_k=u(k/L)$ and 
choose rank one projections 
$p_{k,1},p_{k,2},\dots,p_{k,n}$ so that 
\[ \sum_{i=1}^np_{k,i}=1 \]
and 
\[ u_k=\sum_{i=1}^n\lambda_i(k/L)p_{k,i}. \]
For each $k=0,1,\dots,L-1$, 
there exists a unitary $w_k$ in $M_n$ such that 
\[ w_kp_{k,i}w_k^*=p_{k+1,i} \]
for every $i=1,2,\dots,n$. 
From $\Lip(u)\leq C$ and $\Lip(\lambda_i)\leq C$, we have 
\[ \lVert u_k-u_{k+1}\rVert\leq\frac{C}{L} \]
and 
\[ \lVert w_ku_kw_k^*-u_{k+1}\rVert\leq\frac{C}{L}. \]
It follows that 
\[ \lVert[u_k,w_k]\rVert\leq\frac{2C}{L}<\delta. \]
We claim that $\Bott(u_k,w_k)$ is zero. 
The two paths 
\[ I_k\ni t\mapsto \sum_{i=1}^n\lambda_i(t)p_{k+1,i}\in M_n \]
and $u:I_k\to M_n$ give us a short path $x$ 
from $w_ku_kw_k^*$ to $u_k$. 
It is not so hard to see that 
the rotation number of $\det(x)$ is zero. 
Hence, we obtain $\Bott(u_k,w_k)=0$. 
Therefore, from Lemma \ref{SuperHomotopy}, 
we can find a path of unitaries $z_k:I_k\to M_n$ such that 
\[ z_k(k/L)=1, \quad z_k((k+1)/L)=w_k, \quad \Lip(z_k)\leq4L \]
and 
\[ \lVert[u_k,z_k(t)]\rVert<\ep/2 \]
for all $t\in I_k$. 

We define a path of self-adjoint elements $h_k:I_k\to M_n$ by 
\[ h_k(t)
=z_k(t)\left(\sum_{i=1}^ng_i(t)p_{k,i}\right)z_k^*(t) \]
for $k=0,1,\dots,L-1$ and $t\in I_k$. 
We have 
\[ \Lip(h_k)
\leq\frac{C}{12}\Lip(z_k)+\frac{C}{6}+\frac{C}{12}\Lip(z_k)
\leq C', \]
and 
\begin{align*}
& \lVert u(t)-e^{2\pi\sqrt{-1}h_k(t)}\rVert \\
&\leq\lVert u(t)-u_k\rVert+\lVert u_k-z_k(t)u_kz_k(t)^*\rVert
+\lVert z_k(t)u_kz_k(t)^*-e^{2\pi\sqrt{-1}h_k(t)}\rVert \\
&\leq\frac{C}{L}+\frac{\ep}{2}+\frac{C}{L}<\ep
\end{align*}
for $t\in I_k$. 
By connecting $h_0,h_1,\dots,h_{L-1}$, 
we get a desired path $h:[0,1]\to M_n$. 
\end{proof}

The following is a well-known fact and 
one can find its proof in \cite[Lemma 4.8]{HR} for example. 

\begin{lem}\label{exp}
Let $h_1$ and $h_2$ be self-adjoint elements in a unital $C^*$-algebra. 
Then 
\[ \lVert e^{2\pi\sqrt{-1}h_1}-e^{2\pi\sqrt{-1}h_2}\rVert
\leq2\pi\lVert h_1-h_2\rVert. \]
\end{lem}

The following lemma is also well-known. 
We have been unable to find a suitable reference in the literature, 
so we include a proof for completeness. 

\begin{lem}\label{log}
Let $A$ be a unital $C^*$-algebra and 
let $u_1,u_2$ be unitaries in $A$ 
such that $\lVert u_i-1\rVert<1/2$ for $i=1,2$. 
Put $2\pi\sqrt{-1}h_i=\log u_i$. 
Then we have 
\[ \lVert h_1-h_2\rVert\leq\pi^{-1}\lVert u_1-u_2\rVert. \]
\end{lem}
\begin{proof}
Since $\log z$ is written by the absolutely convergent power series 
\[ \log z=-\sum_{n=1}^\infty\frac{1}{n}(1-z)^n \]
on $\{z\in\C\mid\lvert z-1\rvert\leq1/2\}$, 
one has 
\[ \lVert h_1-h_2\rVert
=\frac{1}{2\pi}\lVert\log u_1-\log u_2\rVert
\leq\frac{1}{2\pi}\sum_{n=1}^\infty\frac{1}{n}
\lVert(1-u_1)^n-(1-u_2)^n\rVert. \]
By an elementary estimate, 
\[ \lVert(1-u_1)^n-(1-u_2)^n\rVert\leq n2^{1-n}\lVert u_1-u_2\rVert \]
is obtained, and so the proof is completed. 
\end{proof}

The following proposition plays a crucial role 
in our cohomology vanishing theorem. 
We let $E=\{t\in\R^2\mid \lVert t\rVert\leq1\}$ and 
$\partial E=\{t\in\R^2\mid \lVert t\rVert=1\}$, 
where $\lVert\cdot\rVert$ is the $\ell^\infty$-norm on $\R^2$. 

\begin{prop}\label{Lipshrink2}
For any $C>0$, 
there exists $C'>0$ such that the following holds. 
Let $A$ be a unital AF algebra and 
let $z:\partial E\to U(A)$ be a continuous map such that 
$z(1,1)=1$ and $\Lip(z)\leq C$. 
If the $K_1$-class of $z$ in $K_1(C(\partial E)\otimes A)$ is zero, 
then we can find a continuous map $\tilde{z}$ 
from $E$ to unitaries of $A$ such that the following are satisfied. 
\begin{enumerate}
\item $z(t)=\tilde{z}(t)$ for all $t\in\partial E$. 
\item $\Lip(\tilde{z})\leq C'$. 
\end{enumerate}
\end{prop}
\begin{proof}
By Proposition \ref{Lipshrink1} and Lemma \ref{exp}, 
there exists $z_0:E\to U(A)$ satisfying the following. 
\begin{itemize}
\item $\lVert z(t)-z_0(t)\rVert<1/2$ for all $t\in\partial E$. 
\item $\Lip(z_0)\leq C_0$, 
where $C_0>0$ is a constant depending only on $C>0$. 
\end{itemize}
Define $k:\partial E\to A$ and $z_1:[0,1]\times\partial E\to U(A)$ by 
\[ k(t)=\frac{1}{2\pi\sqrt{-1}}\log(z(t)^*z_0(t)) \]
and
\[ z_1(s,t)=z(t)e^{2\pi\sqrt{-1}sk(t)} \]
for $t\in\partial E$ and $s\in[0,1]$. 
Clearly $z_1(0,t)=z(t)$ and $z_1(1,t)=z_0(t)$. 
It follows from Lemma \ref{log} and Lemma \ref{exp} that 
$z_1(s,t)$ is Lipschitz continuous 
with some universal Lipschitz constant depending only on $C>0$. 
Define 
\[ \tilde{z}(t)=\begin{cases}
z_1(2-2\lVert t\rVert,t/\lVert t\rVert) & \text{ if }\lVert t\rVert\geq1/2 \\
z_0(2t) & \text{ if }\lVert t\rVert\leq1/2 \end{cases} \]
for $t\in E$. 
Then the assertion easily follows. 
\end{proof}

Now we would like to prove the cohomology vanishing theorem. 

\begin{thm}\label{CVanish1}
Let $A$ be a unital AF algebra and 
let $\alpha$ be an action of $\Z^2$ on $A$ with the Rohlin property. 
Suppose that $\alpha_n$ is approximately inner for any $n\in\Z^2$. 
Then, for any admissible $\alpha$-cocycle $\{u_n\}_n$ in $A^\infty$, 
there exists a unitary $v\in A^\infty$ such that 
\[ u_n=v\alpha_n(v)^* \]
for any $n\in\Z^2$. 
\end{thm}
\begin{proof}
It suffices to show that, for any $\ep>0$, 
there exists a unitary $v\in A^\infty$ such that 
\[ \lVert u_{\xi_i}-v\alpha_{\xi_i}(v)^*\rVert<\ep \]
for each $i=1,2$. 
By applying Proposition \ref{Lipshrink2} to $C=4$, 
we obtain a constant $C'>0$. 
Choose a natural number $M$ so that $2C'/M$ is less than $\ep$. 
Since $\alpha$ has the Rohlin property, 
there exist $R\in\N$ and $m^{(1)},m^{(2)},\dots,m^{(R)}\in\Z^2$ 
with $m^{(1)},\dots,m^{(R)}\geq(M,M)$ 
and which satisfies the requirement in Definition \ref{Rohlin}. 

For each $n\in\Z^2$, 
let $(u_n(k))_{k\in\N}\in\ell^\infty(\N,A)$ be 
a representing sequence of $u_n$. 
We may assume that $u_n(k)$ is a unitary. 
Since $\{u_n\}_n$ is admissible, 
$\kappa(u_{\xi_1}(k),u_{\xi_2}(k),\alpha_{\xi_1},\alpha_{\xi_2})$ 
is zero in $K_0(A)$ for sufficiently large $k\in\N$. 

For each $r=1,2,\dots,R$ and $i=1,2$, 
let $m^{(r)}_i$ be the $i$-th summand of $m^{(r)}$. 
We put $\eta_{r,i}=m^{(r)}_i\xi_i\in\Z^2$. 
For sufficiently large $k\in\N$, 
the pair of unitaries $\{u_{\eta_{r,1}}(k),u_{\eta_{r,2}}(k)\}$ is 
an almost cocycle for $\alpha_{\eta_{r,1}},\alpha_{\eta_{r,2}}$. 
Furthermore, by Lemma \ref{power}, 
$\kappa(u_{\eta_{r,1}}(k),u_{\eta_{r,2}}(k),
\alpha_{\eta_{r,1}},\alpha_{\eta_{r,2}})$ is zero in $K_0(A)$. 
It follows from Lemma \ref{zONboundary1} that 
for sufficiently large $k\in\N$, 
there exists a continuous map $z^{(r)}_k:\partial E\to U(A)$ such that 
the following are satisfied. 
\begin{itemize}
\item $z^{(r)}_k(1,1)=1$. 
\item For all $t\in[-1,1]$, 
\[ \lim_{k\to\infty}\lVert z^{(r)}_k(-1,t)
-u_{\eta_{r,1}}(k)\alpha_{\eta_{r,1}}(z^{(r)}_k(1,t))\rVert=0 \]
and 
\[ \lim_{k\to\infty}\lVert z^{(r)}_k(t,-1)
-u_{\eta_{r,2}}(k)\alpha_{\eta_{r,2}}(z^{(r)}_k(t,1))\rVert=0. \] 
\item $\Lip(z^{(r)}_k)\leq4$. 
\item The $K_1$-class of $z^{(r)}_k$ in $K_1(C(\partial E)\otimes A)$ 
is zero. 
\end{itemize}
From Proposition \ref{Lipshrink2}, 
for sufficiently large $k\in\N$, 
we obtain $\tilde{z}^{(r)}_k:E\to U(A)$ such that 
the following are satisfied. 
\begin{itemize}
\item $z^{(r)}_k(t)=\tilde{z}^{(r)}_k(t)$ for all $t\in\partial E$. 
\item $\Lip(\tilde{z}^{(r)}_k)\leq C'$. 
\end{itemize}
For each $r=1,2,\dots,R$, $g=(g_1,g_2)\in\Z_{m^{(r)}}$ and 
sufficiently large $k\in\N$, 
we define $w^{(r)}_g(k)$ in $U(A)$ by 
\[ w^{(r)}_g(k)=\tilde{z}^{(r)}_k\left(\frac{2g_1}{m^{(r)}_1}-1,
\frac{2g_2}{m^{(r)}_2}-1\right). \]
Let $w^{(r)}_g$ be the image of $(w^{(r)}_g(k))_{k\in\N}$ in $A^\infty$. 
It is easily seen that one has the following 
for any $r=1,2,\dots,R$, $i=1,2$ and $g=(g_1,g_2)\in\Z_{m^{(r)}}$. 
\begin{itemize}
\item If $g_i\neq0$, then 
$\lVert w^{(r)}_g-w^{(r)}_{g-\xi_i}\rVert$ is less than $\ep$. 
\item If $g_i=0$, then 
the distance from $w^{(r)}_g$ to 
$u_{\eta_{r,i}}\alpha_{\eta_{r,i}}(w^{(r)}_{g+\eta_{r,i}-\xi_i})$ 
is less than $\ep$. 
\end{itemize}
We can take a family of projections 
$\{e^{(r)}_g\mid r=1,2,\dots,R, \ g\in\Z_{m^{(r)}}\}$ 
in $A_\infty$ such that 
\[ \sum_{r=1}^R\sum_{g\in\Z_{m^{(r)}}}e^{(r)}_g=1, 
\quad \alpha_{\xi_i}(e^{(r)}_g)=e^{(r)}_{g+\xi_i} \]
for any $r=1,2,\dots,R$, $i=1,2$ and $g\in\Z_{m^{(r)}}$, 
where $g+\xi_i$ is understood modulo $m^{(r)}\Z^2$. 
Moreover, we may also assume that 
$e^{(r)}_g$ commutes with $u_g$ and $\alpha_g(w^{(r)}_g)$. 
Define $v\in\mathcal{U}(A^\infty)$ by 
\[ v=\sum_{r=1}^R\sum_{g\in\Z_{m^{(r)}}}
u_g\alpha_g(w^{(r)}_g)e^{(r)}_g. \]
It is now routinely checked that 
$\lVert u_{\xi_i}-v\alpha_{\xi_i}(v^*)\rVert$ is 
less than $\ep$ for each $i=1,2$. 
\end{proof}

When $A$ is a UHF algebra, 
we can show the following in a similar fashion to the theorem above. 

\begin{thm}\label{CVanish2}
Let $A$ be a UHF algebra and 
let $\alpha$ be an action of $\Z^2$ on $A$ with the Rohlin property. 
Then, for any admissible $\alpha$-cocycle $\{u_n\}_n$ in $A_\infty$, 
there exists a unitary $v\in A_\infty$ such that 
\[ u_n=v\alpha_n(v)^* \]
for any $n\in\Z^2$. 
\end{thm}
\begin{proof}
It suffices to show that, 
for any $\ep>0$ and any unital full matrix subalgebra $A_0\subset A$, 
there exists a unitary $v\in A^\infty\cap A_0'$ such that 
\[ \lVert u_{\xi_i}-v\alpha_{\xi_i}(v)^*\rVert<\ep \]
for each $i=1,2$. 

The proof is almost the same as that of Theorem \ref{CVanish1}. 
Choose $C'>0$, $M\in\N$ and $m^{(1)},m^{(2)},\dots,m^{(R)}\in\Z^2$ 
in the same way. 
By using Lemma \ref{zONboundary2} instead of Lemma \ref{zONboundary1}, 
we may assume that the range of the map $z^{(r)}_k$ is $U(A\cap A_0')$ 
for all $r=1,2,\dots,R$ and sufficiently large $k\in\N$. 
By Lemma \ref{zONboundary2}, 
the $K_1$-class of $z^{(r)}_k$ in $K_1(C(\partial E)\otimes A)$ is zero. 
It follows that 
the $K_1$-class of $z^{(r)}_k$ in $K_1(C(\partial E)\otimes(A\cap A_0'))$ 
is also zero, because $A_0$ is a full matrix subalgebra. 
We apply Proposition \ref{Lipshrink2} to $A\cap A_0'$ and $z^{(r)}_k$, 
and obtain a continuous map $\tilde{z}^{(r)}_k:E\to U(A\cap A_0')$. 
The rest of the proof is exactly the same as Theorem \ref{CVanish1}. 
\end{proof}

The following corollary is an immediate consequence of the theorem above 
and we omit the proof. 
Note that if $A$ is a UHF algebra, then 
any $\alpha$-cocycle in $A$ is admissible. 

\begin{cor}\label{appCVanish}
Let $A$ be a UHF algebra and 
let $\alpha$ be an action of $\Z^2$ on $A$ with the Rohlin property. 
For any $\ep>0$ and a finite subset $\mathcal{F}$ of $A$, 
there exist $\delta>0$ and a finite subset $\mathcal{G}$ of $A$ 
such that the following holds: 
If a family of unitaries $\{u_n\}_{n\in\Z^2}$ in $A$ 
is an $\alpha$-cocycle satisfying 
\[ \lVert[u_{\xi_i},a]\rVert<\delta \]
for each $i=1,2$ and $a\in\mathcal{G}$, then 
we can find a unitary $v\in U(A)$ 
satisfying 
\[ \lVert u_{\xi_i}-v\alpha_{\xi_i}(v^*)\rVert<\ep \]
and 
\[ \lVert[v,a]\rVert<\ep \]
for each $i=1,2$ and $a\in\mathcal{F}$. 
\end{cor}

%%%%%%%%%%%%%%%%%%%%%%%%%%%%%%%%%%%%%%%%%%%%%%%%%%%%%%%%%%%%
\section{Invariants}

In this section, 
we would like to introduce an invariant $[\alpha]$ 
for a $\Z^2$-action $\alpha$ on a UHF algebra. 
Let $A$ be a UHF algebra and 
let $\tau$ be the unique tracial state on $A$. 
Throughout this section, we identify $K_0(A)$ and $\tau(K_0(A))\subset\R$. 
We need the generalized determinant 
introduced by P. de la Harpe and G. Skandalis 
(see \cite{HS} for details). 
For any $u\in U(A)$, 
there exists a piecewise smooth path $h:[0,1]\to U(A)$ 
such that $h(0)=1$ and $h(1)=u$. 
It is not so hard to see that 
\[ \frac{1}{2\pi\sqrt{-1}}\int_0^1\tau(\dot{h}(t)h(t)^*)\,dt \]
is real. 
The determinant $\Delta_\tau$ 
associated with the tracial state $\tau$ is the map 
from $U(A)$ to $\R/\tau(K_0(A))$ defined by 
\[ \Delta_\tau(u)
=\frac{1}{2\pi\sqrt{-1}}\int_0^1\tau(\dot{h}(t)h(t)^*)\,dt
+\tau(K_0(A)). \]
A crucial fact here is that $\Delta_\tau$ is a group homomorphism. 
We also notice that 
$\Delta_\tau$ is $\alpha$-invariant for any $\alpha\in\Aut(A)$, 
that is, $\Delta_\tau\circ\alpha=\Delta_\tau$. 

\bigskip

Let $A$ be a UHF algebra. 
For every prime number $p$, we define 
\[ \zeta(p)=\sup\{k\in\N\cup\{0\}\mid[1]\text{ is divisible by }p^k
\text{ in }K_0(A)\}\in\{0,1,2,\ldots,\infty\}, \]
and let $\mathcal{P}(A)$ be the set of prime numbers $p$ 
such that $1\leq\zeta(p)<\infty$. 
Let $p\in\mathcal{P}(A)$. 
We put $\theta(p)=p^{\zeta(p)}$. 
There exists a unital full matrix subalgebra $A_0\subset A$ 
which is isomorphic to $M_{\theta(p)}$. 
It is clear that 
the canonical map from $K_0(A\cap A_0')$ to $K_0(A)$ is injective 
and $K_0(A)/K_0(A\cap A_0')$ is isomorphic to $\Z/\theta(p)\Z$. 
We also set 
\[ \mathcal{Q}(A)=\{\theta(p)\mid p\in\mathcal{P}(A)\}
\cup\{p^n\mid\text{$p$ with $\zeta(p)=\infty$ and $n\in\N$}\}, \]
so that $A$ is isomorphic to $\bigotimes_{q\in\mathcal{Q}(A)}M_{q}$. 

Let $\alpha,\beta$ be actions of $\Z^2$ on $A$. 
To simplify notation, 
we denote $\alpha_{\xi_i}$, $\beta_{\xi_i}$ 
by $\alpha_i$, $\beta_i$ for each $i=1,2$. 
We would like to introduce an invariant $[\beta,\alpha]$. 
Choose a positive real number $\delta_0>0$ so that, 
for any unital finite dimensional $C^*$-subalgebra $A_0\subset A$ and 
a unitary $w\in A$ satisfying $\lVert[w,a]\rVert<\delta_0$ 
for all $a\in U(A_0)$, 
there exists a piecewise smooth path of unitaries $h:[0,1]\to U(A)$ 
such that $h(0)=1$, $h(1)=w$ and $\lVert[h(t),a]\rVert<1/8$ 
for all $t\in[0,1]$ and $a\in U(A_0)$. 
We may assume $\delta_0<1/8$. 

Take $p\in\mathcal{P}(A)$ arbitrarily. 
Let $A_0\subset A$ be a unital full matrix subalgebra 
which is isomorphic to $M_{\theta(p)}$. 
Take a finite dimensional $C^*$-subalgebra $A_1\subset A$ so that 
for any unitary $w\in U(A_0)$, there exist unitaries $a,b,c\in U(A_1)$ 
such that 
\[ \lVert w-a\rVert<\delta_0/7, \quad 
\lVert \beta_1^{-1}(w)-b\rVert<\delta_0/7 \ \text{ and } \ 
\lVert \beta_2^{-1}(w)-c\rVert<\delta_0/7. \]
We can find unitaries $u_1$ and $u_2$ in $A$ such that 
\[ \lVert\beta_i(a)-\Ad u_i\circ\alpha_i(a)\rVert<\delta_0/7 \]
for all $a$ in $\beta_1^{-1}(U(A_1))\cup\beta_2^{-1}(U(A_1))
\cup\beta_1^{-1}\beta_2^{-1}(U(A_1))$. 
Define $x\in U(A)$ by 
\[ x=u_1\alpha_1(u_2)(u_2\alpha_2(u_1))^*. \]
It is easy to see that 
$\lVert[x,a]\rVert$ is less than $\delta_0$ for every $a\in U(A_0)$. 
There exists a piecewise smooth path of unitaries $h:[0,1]\to U(A)$ 
such that $h(0)=1$, $h(1)=x$ and 
$\lVert[h(t),a]\rVert<1/2$ for all $t\in[0,1]$ and $a\in U(A_0)$. 
Since $\Delta_\tau$ is a group homomorphism and 
$\alpha_i$-invariant for each $i=1,2$, one has $\Delta_\tau(x)=0$. 
Hence, 
\[ \frac{1}{2\pi\sqrt{-1}}\int_0^1\tau(\dot{h}(t)h(t)^*)\,dt
\in K_0(A). \]
We denote the image of this value 
in $K_0(A)/K_0(A\cap A_0')\cong\Z/\theta(p)\Z$ by $[\beta,\alpha](p)$. 
Thus, 
\begin{align*}
[\beta,\alpha](p)
&=\frac{1}{2\pi\sqrt{-1}}\int_0^1\tau(\dot{h}(t)h(t)^*)\,dt
+K_0(A\cap A_0') \\ 
&\in K_0(A)/K_0(A\cap A_0')\cong\Z/\theta(p)\Z. 
\end{align*}

\begin{lem}
In the setting above, 
$[\beta,\alpha](p)$ does not depend on the choice of 
$A_0$, $A_1$, $u_1$, $u_2$ and $h$. 
\end{lem}
\begin{proof}
Let $k:[0,1]\to U(A)$ be another piecewise smooth path such that 
$k(0)=1$, $k(1)=x$ and 
$\lVert[k(t),a]\rVert<1/2$ for all $t\in[0,1]$ and $a\in U(A_0)$. 
By connecting $h$ and $k$, 
we get a closed path $H$ which almost commutes with elements in $U(A_0)$. 
Clearly $H$ is homotopic to a closed path in $U(A\cap A_0')$, 
and so the difference between 
\[ \frac{1}{2\pi\sqrt{-1}}\int_0^1\tau(\dot{h}(t)h(t)^*)\,dt \]
and 
\[ \frac{1}{2\pi\sqrt{-1}}\int_0^1\tau(\dot{k}(t)k(t)^*)\,dt \]
is in $K_0(A\cap A_0')$. 
Therefore, $[\beta,\alpha](p)\in K_0(A)/K_0(A\cap A_0')$ 
does not depend on the choice of the path $h$. 

Next, we would like to verify that 
$[\beta,\alpha](p)$ does not depend on the choice of $u_1$. 
Suppose that $u_1'\in U(A)$ also satisfies 
\[ \lVert\beta_1(a)-\Ad u_1'\circ\alpha_1(a)\rVert<\delta_0/7 \]
for all $a$ in $\beta_1^{-1}(U(A_1))\cup\beta_2^{-1}(U(A_1))
\cup\beta_1^{-1}\beta_2^{-1}(U(A_1))$. 
Put $v=u_1'u_1^*$. 
Then we have $\lVert[a,v]\rVert<2\delta_0/7$ for all $a\in U(A_1)$. 
Hence we can find a piecewise smooth path of unitaries $k:[0,1]\to U(A)$ 
such that $k(0)=1$, $k(1)=v$ and 
$\lVert[k(t),a]\rVert<1/8$ for all $t\in[0,1]$ and $a\in U(A_1)$. 
It follows that 
\[ \lVert[k(t),a]\rVert<\frac{2\delta_0}{7}+\frac{1}{8} \ \text{ and } \ 
\lVert[k(t),\beta_2^{-1}(a)]\rVert<\frac{2\delta_0}{7}+\frac{1}{8} \]
for all $t\in[0,1]$ and $a\in U(A_0)$. 
For the unitary $x=u_1\alpha_1(u_2)(u_2\alpha_2(u_1))^*$, 
there exists a piecewise smooth path of unitaries $h:[0,1]\to U(A)$ 
such that $h(0)=1$, $h(1)=x$ and 
$\lVert[h(t),a]\rVert<1/8$ for all $t\in[0,1]$ and $a\in U(A_0)$. 
Define $H:[0,1]\to U(A)$ by 
\[ H(t)=k(t)h(t)u_2\alpha_2(k(t)^*)u_2^*. \]
We can see that 
\[ H(0)=1, \ H(1)=u_1'\alpha_1(u_2)(u_2\alpha_2(u_1'))^* \]
and 
\[ \lVert[H(t),a]\rVert<\frac{2\delta_0}{7}+\frac{1}{8}+\frac{1}{8}
+\frac{5\delta_0}{7}+\frac{1}{8}=\delta_0+3/8<1/2 \]
for all $t\in[0,1]$ and $a\in U(A_0)$. 
In addition, by \cite[Lemma 1]{HS}, 
\[ \frac{1}{2\pi\sqrt{-1}}\int_0^1\tau(\dot{H}(t)H(t)^*)\,dt
=\frac{1}{2\pi\sqrt{-1}}\int_0^1\tau(\dot{h}(t)h(t)^*)\,dt \]
Therefore, $[\beta,\alpha](p)$ does not depend on the choice of $u_1$. 
Likewise, we can show that 
$[\beta,\alpha](p)$ does not depend on the choice of $u_2$. 

In order to show that 
$[\beta,\alpha](p)$ does not depend on the choice of $A_1$, 
we take another finite dimensional $C^*$-subalgebra $A_2$ so that, 
for any unitary $w\in U(A_0)$, there exist unitaries $a,b,c\in U(A_2)$ 
such that 
\[ \lVert w-a\rVert<\delta_0/7, \quad 
\lVert \beta_1^{-1}(w)-b\rVert<\delta_0/7 \ \text{ and } \ 
\lVert \beta_2^{-1}(w)-c\rVert<\delta_0/7. \]
Since we can find unitaries $u_1,u_2$ in $A$ such that 
\[ \lVert\beta_i(a)-\Ad u_i\circ\alpha_i(a)\rVert<\delta_0/7 \]
for all $a$ in 
\[ \beta_1^{-1}(U(A_1)\cup U(A_2))\cup\beta_2^{-1}(U(A_1)\cup U(A_2))
\cup\beta_1^{-1}\beta_2^{-1}(U(A_1)\cup U(A_2)), \]
the assertion is clear. 

Finally, we observe that 
$[\beta,\alpha](p)$ does not depend on the choice of 
$A_0\cong M_{\theta(p)}$. 
Suppose that $B_0$ is also a full matrix subalgebra of $A$ 
which is isomorphic to $M_{\theta(p)}$. 
By taking the finite dimensional $C^*$-subalgebra $A_1$ so large, 
we may assume that 
$x=u_1\alpha_1(u_2)(u_2\alpha_2(u_1))^*$ can be connected to the identity 
by a piecewise smooth path of unitaries 
which almost commute with both $A_0$ and $B_0$. 
Hence, $[\beta,\alpha](p)$ does not depend on the choice of $A_0$. 
\end{proof}

We put $[\beta,\alpha]\in\prod_{p\in\mathcal{P}(A)}\Z/\theta(p)\Z$ by 
\[ [\beta,\alpha]=([\beta,\alpha](p))_{p\in\mathcal{P}(A)}. \]
Furthermore, we set $[\alpha]=[\id,\alpha]$. 
We will show that 
$[\alpha]$ is the complete invariant for cocycle conjugacy 
in the next section. 

\begin{lem}\label{additive}
Let $A$ be a UHF algebra and 
let $\alpha,\beta,\gamma$ be actions of $\Z^2$ on $A$. 
Then one has $[\gamma,\alpha]=[\gamma,\beta]+[\beta,\alpha]$. 
\end{lem}
\begin{proof}
Take $p\in\mathcal{P}(A)$ arbitrarily. 
We will prove $[\gamma,\alpha](p)=[\gamma,\beta](p)+[\beta,\alpha](p)$. 
To simplify notation, 
we denote $\alpha_{\xi_i}$, $\beta_{\xi_i}$, $\gamma_{\xi_i}$ 
by $\alpha_i$, $\beta_i$, $\gamma_i$ for each $i=1,2$. 
Let $A_0\subset A$ be a unital full matrix subalgebra 
which is isomorphic to $M_{\theta(p)}$. 
Choose a sufficiently large finite subset $\mathcal{F}\subset A$ 
and a sufficiently small positive real number $\delta>0$. 
There exist unitaries $v_1,v_2$ such that 
\[ \lVert\gamma_i(a)-\Ad v_i\circ\beta_i(a)\rVert<\delta \]
for all $a\in\mathcal{F}$ and $i=1,2$. 
We may also find piecewise smooth paths of unitaries 
$u_1,u_2:[1,\infty)\to U(A)$ such that 
\[ \beta_i(a)=\lim_{s\to\infty}\Ad u_i(s)\circ\alpha_i(a) \]
for all $a\in A$. 
We may further assume 
\[ \lVert\beta_i(a)-\Ad u_i(s)\circ\alpha_i(a)\rVert<\delta \]
for all $a\in\mathcal{F}\cup\{v_1,v_2\}$, $s\in[1,\infty)$ and $i=1,2$. 
We write $u_1=u_1(1)$ and $u_2=u_2(1)$. 
Define $x,y\in U(A)$ by 
\[ x=u_1\alpha_1(u_2)(u_2\alpha_2(u_1))^* \]
and 
\[ y=v_1\beta_1(v_2)(v_2\beta_2(v_1))^*. \]
There exist piecewise smooth paths of unitaries 
$h,k:[0,1]\to U(A)$ such that 
\[ h(0)=k(0)=1, \ h(1)=x, \ k(1)=y \]
and $h(t),k(t)$ almost commute with unitaries in $A_0$. 
We have 
\begin{align*}
&[\gamma,\beta](p)+[\beta,\alpha](p) \\
&=\frac{1}{2\pi\sqrt{-1}}\int_0^1\tau(\dot{h}(t)h(t)^*)\,dt
+\frac{1}{2\pi\sqrt{-1}}\int_0^1\tau(\dot{k}(t)k(t)^*)\,dt
+K_0(A\cap A_0')
\end{align*}

Since 
\[ \lVert\gamma_i(a)-\Ad v_iu_i\circ\alpha_i(a)\rVert<2\delta \]
for all $a\in\mathcal{F}$ and $i=1,2$, 
\[ z=v_1u_1\alpha_1(v_2u_2)(v_2u_2\alpha_2(v_1u_1))^* \]
almost commutes with unitaries in $A_0$. 
We may also assume that 
\[ z(t)=(v_1u_1\alpha_1(v_2)u_1^*)h(t)(v_2u_2\alpha_2(v_1)u_2^*)^* \]
almost commutes with unitaries in $A_0$ for every $t\in[0,1]$. 
It is straightforward to see $z(1)=z$. 
Let $w:[0,1]\to U(A)$ be a path of unitaries such that 
$w(0)=0$, $w(1)=z(0)$ and 
$w(t)$ almost commutes with unitaries in $A_0$. 
Then, one has 
\begin{align*}
&[\gamma,\alpha](p) \\
&=\frac{1}{2\pi\sqrt{-1}}\int_0^1\tau(\dot{w}(t)w(t)^*)\,dt
+\frac{1}{2\pi\sqrt{-1}}\int_0^1\tau(\dot{z}(t)z(t)^*)\,dt
+K_0(A\cap A_0') \\
&=\frac{1}{2\pi\sqrt{-1}}\int_0^1\tau(\dot{w}(t)w(t)^*)\,dt
+\frac{1}{2\pi\sqrt{-1}}\int_0^1\tau(\dot{h}(t)h(t)^*)\,dt
+K_0(A\cap A_0'). 
\end{align*}
Therefore, it suffices to show 
\[ \frac{1}{2\pi\sqrt{-1}}\left(\int_0^1\tau(\dot{k}(t)k(t)^*)\,dt
-\int_0^1\tau(\dot{w}(t)w(t)^*)\,dt\right)
\in K_0(A\cap A_0'). \]

By connecting $k$, $w$ and a short path from $y$ to $z(0)$, 
we obtain a closed path $H$ 
which is almost commuting with unitaries in $A_0$. 
Let $\lambda\in K_0(A\cap A_0')\subset\R$ be 
the $K_0$-class determined by the closed path $H$. 
We would like to show 
\[ \lambda=\frac{1}{2\pi\sqrt{-1}}\left(\int_0^1\tau(\dot{k}(t)k(t)^*)\,dt
-\int_0^1\tau(\dot{w}(t)w(t)^*)\,dt\right). \]
For $s\in[1,\infty)$, we let 
\[ f(s)=(v_1u_1(s)\alpha_1(v_2)u_1(s)^*)(v_2u_2(s)\alpha_2(v_1)u_2(s)^*)^* \]
Thus, $f$ is a path from $z(0)$ to $y$. 
For any $\ep>0$, 
there exists $s\in[1,\infty)$ such that 
\[ \lVert\beta_1(v_2)-\Ad u_1(s)\circ\alpha_1(v_2)\rVert<\ep \]
and 
\[ \lVert\beta_2(v_1)-\Ad u_2(s)\circ\alpha_2(v_1)\rVert<\ep. \]
We define 
\[ a=\frac{1}{2\pi\sqrt{-1}}\log(f(s)y^*). \]
Then $\lVert a\rVert$ is less than $\arcsin2\ep$. 
By connecting paths $k$, $w$, $f:[1,s]\to U(A)$ and 
\[ g:t\mapsto e^{2\pi\sqrt{-1}ta}y, \]
we obtain a closed path $H'$. 
Clearly $H$ is homotopic to $H'$, because they are close to each other. 
Moreover, we can see 
\[ \frac{1}{2\pi\sqrt{-1}}\int_1^s\tau(\dot{f}(t)f(t)^*)\,dt=0 \]
and 
\[ \left\lvert\frac{1}{2\pi\sqrt{-1}}\int_0^1
\tau(\dot{g}(t)g(t)^*)\,dt\right\rvert=\lvert\tau(a)\rvert<\arcsin2\ep. \]
It follows that 
\[ \left\lvert\lambda
-\frac{1}{2\pi\sqrt{-1}}\left(\int_0^1\tau(\dot{k}(t)k(t)^*)\,dt
-\int_0^1\tau(\dot{w}(t)w(t)^*)\,dt\right)\right\rvert<\arcsin2\ep. \]
Since $\ep>0$ was arbitrary, we get the conclusion. 
\end{proof}

We need to recall the notion of outer conjugacy. 
Two actions $\alpha$ and $\beta$ of $\Z^2$ on $A$ is 
said to be outer conjugate, 
if there exist an automorphism $\sigma\in\Aut(A)$ and 
unitaries $\{u_n\}_{n\in\Z^2}\subset U(A)$ such that 
\[ \Ad u_n\circ\alpha_n(a)=\sigma\circ\beta_n\circ\sigma^{-1}(a) \]
for any $n\in\Z^2$ and $a\in A$. 
Notice that $\{u_n\}_n$ is not necessarily an $\alpha$-cocycle. 

The following lemma is essentially contained in \cite{N1}. 

\begin{lem}\label{model}
Let $A$ be a UHF algebra. 
For any $f\in\prod_{p\in\mathcal{P}(A)}\Z/\theta(p)\Z$, 
we can construct a uniformly outer $\Z^2$-action $\gamma^f$ on $A$ 
such that $[\gamma^f]=f$. 
Moreover, if $f,g\in\prod\Z/\theta(p)\Z$ satisfy 
$f(p)=g(p)$ for all but finitely many $p\in\mathcal{P}(A)$, then 
$\gamma^f$ is outer conjugate to $\gamma^g$. 
\end{lem}
\begin{proof}
Let 
\[ A=\bigotimes_{q\in\mathcal{Q}(A)}M_q \]
be a tensor product decomposition. 
For $q\in\mathcal{Q}(A)$, 
we define unitaries $u_q,v_q\in M_{q}$ by 
\[ u_q=\begin{bmatrix}
1 & & & & \\
 & \omega & & & \\
 & & \omega^2 & & \\
 & & & \ddots & \\
 & & & & \omega^{q-1}
\end{bmatrix},
\qquad
v_q=\begin{bmatrix}
0 & \cdot & \cdots & \cdot & 0 & 1 \\
1 & 0 & \cdots & & \cdot & 0 \\
0 & 1 & \cdot & & & \cdot \\
\vdots & & \ddots & \ddots & \vdots & \vdots \\
\cdot & & 0 & 1 & 0 & 0 \\
0 & \cdot & \cdots & 0 & 1 & 0
\end{bmatrix}, \]
where $\omega=e^{2\pi\sqrt{-1}/q}$. 
Take a decomposition $\mathcal{Q}(A)=L_1\cup L_2$ 
such that both $L_1$ and $L_2$ are infinite sets. 

Take $f\in\prod_{p\in\mathcal{P}(A)}\Z/\theta(p)\Z$. 
For each $p\in \mathcal{P}(A)$, 
we regard $f(p)$ as an element of $\{0,1,\dots,\theta(p)-1\}$. 
Define a $\Z^2$-action $\gamma^f$ on $A$ by 
\[ \gamma^f_{\xi_1}
=\!\bigotimes_{q=\theta(p)\in L_f}\!\Ad u_q^{f(p)}
\otimes\bigotimes_{q\in L_1\setminus L_f}\id
\otimes\bigotimes_{q\in L_2\setminus L_f}\Ad u_q. \]
and 
\[ \gamma^f_{\xi_2}
=\bigotimes_{q\in L_f}\Ad v_q
\otimes\bigotimes_{q\in L_1\setminus L_f}\Ad v_q
\otimes\bigotimes_{q\in L_2\setminus L_f}\id \]
where $L_f=\{\theta(p)\in \mathcal{Q}(A)\mid f(p) \neq 0\}$. 
We can verify that 
$\gamma^f$ has the Rohlin property and $[\gamma^f]=f$ 
(see the proof of \cite[Theorem~15~(1)]{N1}). 
The latter assertion is easy to see. 
\end{proof}

%%%%%%%%%%%%%%%%%%%%%%%%%%%%%%%%%%%%%%%%%%%%%%%%%%%%%%%%%%%%
\section{Classification}

Let $A$ be a UHF algebra and 
let $\tau$ be the unique tracial state on $A$. 
Throughout this section, we identify $K_0(A)$ and $\tau(K_0(A))\subset\R$. 
We begin with the following lemma. 

\begin{lem}\label{stability}
Let $A$ be a UHF algebra and 
let $\alpha\in\Aut(A)$ be an automorphism with the Rohlin property. 
Then, for any unitary $u\in A_\infty$, 
there exists a unitary $v\in A_\infty$ such that 
\[ u=v\alpha(v^*). \]
\end{lem}
\begin{proof}
This is the cohomology vanishing for single automorphisms. 
The reader should see \cite{HO} and \cite[Lemma 2.9]{I1}. 
\end{proof}

By means of the lemma above, one obtains the following. 

\begin{lem}\label{existing3}
Let $A$ be a UHF algebra and 
let $\alpha,\beta$ be actions of $\Z^2$ on $A$ with the Rohlin property. 
Suppose that there exists a unitary $u_1\in A$ such that 
$\beta_{\xi_1}(a)=\Ad u_1\circ\alpha_{\xi_1}(a)$ for all $a\in A$. 
If $[\alpha]=[\beta]$, then 
for any finite subset $\mathcal{F}\subset A$ and $\ep>0$, 
there exists $u_2\in U(A)$ satisfying the following. 
\begin{enumerate}
\item $\lVert\beta_{\xi_2}(a)-\Ad u_2\circ\alpha_{\xi_2}(a)\rVert$ 
is less than $\ep$ for all $a\in\mathcal{F}$. 
\item $\lVert u_1\alpha_{\xi_1}(u_2)-u_2\alpha_{\xi_2}(u_1)\rVert$ 
is less than $\ep$. 
\item $\kappa(u_1,u_2,\alpha_{\xi_1},\alpha_{\xi_2})$ is zero. 
\end{enumerate}
\end{lem}
\begin{proof}
To simplify notation, 
we denote $\alpha_{\xi_i}$, $\beta_{\xi_i}$ 
by $\alpha_i$, $\beta_i$ for each $i=1,2$. 
First, we claim the following: 
there exists a unitary $u_2\in A^\infty$ such that 
$u_1\alpha_1(u_2)=u_2\alpha_2(u_1)$ and 
$\beta_2(a)=\Ad u_2\circ\alpha_2(a)$ for all $a\in A$. 
Indeed, we can find $x\in U(A^\infty)$ such that 
$\beta_2(a)=\Ad x\circ\alpha_2(a)$ for all $a\in A$. 
It is easy to see that 
$y=u_1\alpha_1(x)(x\alpha_2(u_1))^*$ is a unitary in $A_\infty$. 
Hence, from Lemma \ref{stability}, 
one can find a unitary $z\in A_\infty$ such that 
\[ y^*=z(\Ad u_1\circ\alpha_1)(z^*). \]
Then, $u_2=z^*x\in A^\infty$ meets the requirements. 

Choose $\delta>0$ so that $4\delta+\arcsin\delta<\ep$. 
Let 
\[ A=\bigotimes_{q\in\mathcal{Q}(A)}M_q \]
be a tensor product decomposition. 
Let $\omega\in\Aut(A)$ be a product type automorphism such that 
\[ \omega(e^{(q)}_{i,j})=e^{(q)}_{i+1,j+1}, \]
where $(e^{(q)}_{i,j})_{i,j}$ is a system of matrix units in $M_q$. 
Then, $\omega$ has the Rohlin property. 
Hence, by Theorem \ref{cc}, 
there exists a unitary $w$ and $\sigma\in\Aut(A)$ such that 
\[ \lVert w-1\rVert<\delta \ \text{ and } \ 
\Ad wu_1\circ\alpha_1=\sigma\circ\omega\circ\sigma^{-1}. \]
Choose a full matrix subalgebra $A_0\cong M_r\subset A$ 
so that the following hold. 
\begin{itemize}
\item $A_0$ is a tensor product of finitely many $M_q$'s above. 
\item For any $a\in\mathcal{F}$, there exists $b\in A_0$ 
such that $\lVert\sigma^{-1}(\beta_2(a))-b\rVert<\delta$. 
\end{itemize}

Let $\delta'>0$ be a small positive number and 
let $\mathcal{G}\subset A$ be a finite subset containing $\mathcal{F}$. 
From the claim above, 
we can find a unitary $u_2\in A$ such that 
\begin{itemize}
\item $\lVert\beta_2(a)-\Ad u_2\circ\alpha_2(a)\rVert$ 
is less than $\delta'$ for all $a\in\mathcal{G}$. 
\item $\lVert u_1\alpha_1(u_2)-u_2\alpha_2(u_1)\rVert$ 
is less than $\delta$. 
\end{itemize}
Put 
\[ a=\frac{1}{2\pi\sqrt{-1}}\log(u_1\alpha_1(u_2)(u_2\alpha_2(u_1))^*). \]
We let $\kappa(u_1,u_2,\alpha_1,\alpha_2)=m/l$, 
where $m\in\Z$ and $l\in\N$ are relatively prime. 
By Lemma \ref{tracevalue}, 
$\lvert m/l\rvert=\lvert\tau(a)\rvert$ 
is less than $(2\pi)^{-1}\arcsin\delta$. 

Let $P=\{p\in\mathcal{P}(A)\mid p\text{ divides }r\}$. 
For $p\in P$, let $A_p\subset A_0$ be a unital full matrix subalgebra 
which is isomorphic to $M_{\theta(p)}$. 
By virtue of Lemma \ref{additive}, we have $[\beta,\alpha]=0$. 
It follows that, 
by taking $\mathcal{G}$ so large and $\delta'>0$ so small, 
we may assume that 
\[ [\beta,\alpha](p)=\tau(a)+K_0(A\cap A_p') \]
for all $p\in P$. 
Hence $\tau(a)$ belongs to $K_0(A\cap A_p')$ for each $p\in P$. 
Thus, $\tau(a)=m/l$ is in $K_0(A\cap A_0')$, and so 
there exists a projection $e$ in $A\cap A_0'$ such that 
\[ e=\omega^l(e), \ e+\omega(e)+\dots+\omega^{l-1}(e)=1. \]
Define a unitary $v\in A$ by 
\[ v=\sum_{j=0}^{l-1}e^{2\pi\sqrt{-1}jm/l}\sigma(\omega^j(e)). \]
It is easily seen that 
\[ \lVert v-\Ad wu_1\circ\alpha_1(v)\rVert
<2\pi\lvert m/l\rvert<\arcsin\delta \]
and 
\[ \lVert[\beta_2(a),v]\rVert<2\delta \]
for all $a\in\mathcal{F}$. 
Besides, it can be also checked that 
\[ \kappa(1,v,\Ad wu_1\circ\alpha_1,\id)=-m/l. \]

We put $u_1'=wu_1$ and $u_2'=vu_2$. 
Clearly we have 
\[ \lVert\beta_i(a)-u_i'\alpha_i(a)u_i'^*\rVert<4\delta \]
for all $a\in\mathcal{F}$ and $i=1,2$. 
In addition, 
\[ \lVert u_1'\alpha_1(u_2')-u_2'\alpha_2(u_1')\rVert
<4\delta+\arcsin\delta. \]
From Lemma \ref{adjust}, one also gets 
\begin{align*}
\kappa(u_1',u_2',\alpha_1,\alpha_2)
&=\kappa(wu_1,vu_2,\alpha_1,\alpha_2) \\
&=\kappa(wu_1,u_2,\alpha_1,\alpha_2)+\kappa(1,v,\Ad wu_1\circ\alpha_1,\id) \\
&=\kappa(u_1,u_2,\alpha_1,\alpha_2)-m/l=0, 
\end{align*}
thereby completing the proof. 
\end{proof}

\begin{lem}
Let $A$ be a UHF algebra and 
let $\alpha,\beta$ be actions of $\Z^2$ on $A$ with the Rohlin property. 
If $[\alpha]=[\beta]$, then 
there exists an admissible $\alpha$-cocycle 
$\{u_n\}_{n\in\Z^2}$ in $A^\infty$ such that 
\[ \beta_n(a)=\Ad u_n\circ\alpha_n(a) \]
for all $n\in\Z^2$ and $a\in A$. 
\end{lem}
\begin{proof}
To simplify notation, 
we denote $\alpha_{\xi_i}$, $\beta_{\xi_i}$ by $\alpha_i$, $\beta_i$. 
It suffices to show the following. 
For any $\ep>0$ and a finite subset $\mathcal{F}$ of $A$, 
there exist unitaries $u_1,u_2$ in $A$ 
such that the following are satisfied. 
\begin{itemize}
\item $\lVert\beta_i(a)-\Ad u_i\circ\alpha_i(a)\rVert$ 
is less than $\ep$ for all $a\in\mathcal{F}$ and $i=1,2$. 
\item $\lVert u_1\alpha_1(u_2)-u_2\alpha_2(u_1)\rVert$ is less than $\ep$. 
\item $\kappa(u_1,u_2,\alpha_1,\alpha_2)$ is zero. 
\end{itemize}

By Theorem \ref{cc}, 
there exists a unitary $w_1\in A$ and $\sigma\in\Aut(A)$ such that 
\[ \beta_1=\Ad w_1\circ\sigma\circ\alpha_1\circ\sigma^{-1}. \]
Let $\alpha_i'=\sigma\circ\alpha_i\circ\sigma^{-1}$ for $i=1,2$. 
By Lemma \ref{existing3}, 
there exists $w_2\in U(A)$ such that 
\begin{itemize}
\item $\lVert\beta_2(a)-\Ad w_2\circ\alpha_2'(a)\rVert$ 
is less than $\ep/2$ for all $a\in\mathcal{F}$. 
\item $\lVert w_1\alpha_1'(w_2)-w_2\alpha_2'(w_1)\rVert$ 
is less than $\ep/2$. 
\item $\kappa(w_1,w_2,\alpha_1',\alpha_2')$ is zero. 
\end{itemize}
We can find a unitary $v\in U(A)$ such that 
\[ \lVert\sigma(a)-vav^*\rVert<\ep/8 \]
for all $a$ in 
\[ \sigma^{-1}(\mathcal{G})\cup\alpha_1(\sigma^{-1}(\mathcal{G}))
\cup\alpha_2(\sigma^{-1}(\mathcal{G})), \]
where $\mathcal{G}=\mathcal{F}\cup\{w_1,w_2\}$. 
Then 
the unitaries $u_1=w_1v\alpha_1(v^*)$ and $u_2=w_2v\alpha_2(v^*)$ satisfy 
\[ \lVert\beta_i(a)-\Ad u_i\circ\alpha_i(a)\rVert<\ep \]
for all $a\in\mathcal{F}$ and $i=1,2$, and 
\[ \lVert u_1\alpha_1(u_2)-u_2\alpha_2(u_1)\rVert<\ep. \]
Besides, one can see 
\begin{align*}
\kappa(u_1,u_2,\alpha_1,\alpha_2)
&=\kappa(w_1v\alpha_1(v^*),w_2v\alpha_2(v^*),\alpha_1,\alpha_2) \\
&=\kappa(v^*w_1v,v^*w_2v,\alpha_1,\alpha_2) \\
&=\kappa(\sigma(v^*w_1v),\sigma(v^*w_2v),\alpha'_1,\alpha'_2) \\
&=\kappa(w_1,w_2,\alpha'_1,\alpha'_2) \\
&=0. 
\end{align*}
This completes the proof. 
\end{proof}

Combining the lemma above and Theorem \ref{CVanish1}, 
we can prove the following. 

\begin{prop}\label{generalUHF}
Let $A$ be a UHF algebra and 
let $\alpha,\beta$ be actions of $\Z^2$ on $A$ with the Rohlin property. 
If $[\alpha]=[\beta]$, then 
for any $\ep>0$ and a finite subset $\mathcal{F}$ of $A$, 
there exists an $\alpha$-cocycle $\{u_n\}_{n\in\Z^2}$ 
(in fact a coboundary) in $A$ such that 
\[ \lVert\beta_{\xi_i}(a)-\Ad u_{\xi_i}\circ\alpha_{\xi_i}(a)\rVert
<\ep \]
for all $a\in\mathcal{F}$ and $i=1,2$. 
\end{prop}

\bigskip

Now we are ready to prove our main results. 
The key ingredient of the proof is 
the Evans-Kishimoto intertwining argument \cite[Theorem 4.1]{EK}. 
See also \cite[Theorem 5]{N2}, \cite[Theorem 3.5]{I2} 
and \cite[Theorem 5.2]{M} for this argument. 

\begin{thm}\label{main1}
Let $A$ be a UHF algebra and 
let $\alpha$ and $\beta$ be two uniformly outer actions of $\Z^2$ on $A$. 
Then, the following conditions are equivalent. 
\begin{enumerate}
\item $[\alpha]=[\beta]$. 
\item $\alpha$ is cocycle conjugate to $\beta$. 
\end{enumerate}
\end{thm}
\begin{proof}
The implication from (2) to (1) is clear. 
It suffices to show that (1) implies (2). 
Note that, by Theorem \ref{RohlinType}, 
both $\alpha$ and $\beta$ have the Rohlin property. 
The proof is carried out by the Evans-Kishimoto intertwining argument. 
In each step of this argument, 
we have to use Proposition \ref{generalUHF} and 
Corollary \ref{appCVanish} repeatedly. 
We omit the detail, 
because it is exactly the same as \cite[Theorem 5.2]{M}. 
\end{proof}

A UHF algebra $A$ is said to be of infinite type, 
if $A$ is isomorphic to $A\otimes A$. 
This is equivalent to say that $\mathcal{P}(A)$ is empty. 
As an immediate consequence, we get the following corollary. 

\begin{cor}
Let $A$ be a UHF algebra of infinite type. 
Then, any two uniformly outer actions of $\Z^2$ on $A$ 
are cocycle conjugate. 
\end{cor}

\begin{rem}
From Theorem \ref{main1} and Theorem \ref{CVanish1}, 
one can actually show the following: 
Let $\alpha$ and $\beta$ be uniformly outer actions of $\Z^2$ 
on a UHF algebra $A$. 
If $[\alpha]=[\beta]$, then 
for any $\ep>0$, there exist $\sigma\in\Aut(A)$ and 
an $\alpha$-cocycle $\{u_n\}_{n\in\Z^2}$ such that 
\[ \Ad u_n\circ\alpha_n(a)=\sigma\circ\beta_n\circ\sigma^{-1}(a), \]
\[ \lVert u_{\xi_i}-1\rVert<\ep \]
for any $n\in\Z^2$, $a\in A$ and $i=1,2$. 
\end{rem}

Finally, we would like to give a classification result 
up to outer conjugacy. 

\begin{thm}\label{main2}
Let $A$ be a UHF algebra and 
let $\alpha$ and $\beta$ be two uniformly outer actions of $\Z^2$ on $A$. 
Then, the following conditions are equivalent. 
\begin{enumerate}
\item $[\alpha](p)=[\beta](p)$ for all but finitely many $p\in\mathcal{P}(A)$. 
\item $\alpha$ is outer conjugate to $\beta$. 
\end{enumerate}
\end{thm}
\begin{proof}
This follows from Theorem \ref{main1} and Lemma \ref{model}. 
\end{proof}

\end{document}